\documentclass[12pt, leqno]{amsart}
\usepackage[mathscr]{eucal}
\usepackage{amssymb,amsmath}
\usepackage{amsfonts}
\usepackage{a4}

\newtheorem{theorem}{Theorem}

\newtheorem{corollary}{Corollary}

\begin{document}
\title{Maximal commutative subalgebras, Poisson geometry and Hochschild homology.}

\author[Tomasz Maszczyk]{Tomasz Maszczyk\dag}
\address{Institute of Mathematics\\
Polish Academy of Sciences\\
Sniadeckich 8\newline 00--956 Warszawa, Poland\\
\newline Institute of Mathematics\\
University of Warsaw\\ Banacha 2\newline 02--097 Warszawa, Poland}
\email{t.maszczyk@uw.edu.pl}

\thanks{\dag The author was partially supported by KBN grants:\\
1P03A 036 26 and  115/E-343/SPB/6.PR UE/DIE 50/2005-2008.}

\begin{abstract} A Poisson geometry arising from maximal
commutative subalgebras is studied. A spectral sequence convergent
to Hochschild homology with coefficients in a bimodule is
presented. It depends on the choice of a maximal commutative
subalgebra inducing appropriate filtrations. Its $E^{2}_{p,
q}$-groups are computed in terms of canonical homology with values
in a Poisson module defined by a given bimodule and a maximal
commutative subalgebra.
\end{abstract}

\maketitle

\paragraph{\textbf{1. Introduction}} In the paper \cite{Bry}
Jean-Luc Brylinski introduces a spectral sequence convergent to
Hochschild homology of an almost commutative algebra (filtered
algebra whose associated graded  algebra is commutative) and
computes (in the case when the associated graded  algebra is
smooth) its $E^{2}_{p,q}$-groups in terms of canonical homology of
the associated graded  algebra. As the main applications serve
there envelopping algebras of Lie algebras and algebras of
differential operators of commutative algebras.

In this paper we show another general construction of almost
commutative algebras. They arise always as an effect of choosing a
maximal commutative subalgebra in an arbitrary associative
algebra. We study the Poisson geometry related to this
construction. It turns out that this geometry describes a
nonlinear involutive distribution on the spectrum of the maximal
commutative subalgebra. Every bimodule over the associative
algebra defines a graded sheaf with a flat connection along this
nonlinear involutive distribution. After a slight generalization
of the result of Brylinski we construct a spectral sequence
convergent to Hochschild homology of this almost commutative
algebra with coefficients in an almost symmetric bimodule. We
compute its $E^{2}_{p,q}$-groups in terms of canonical homology of
the associated graded Poisson algebra with values in an associated
graded Poisson module.

The canonical complex was investigated by Gelfand-Dorfman
\cite{Gel}, Koszul \cite{Kos}, Brylinski \cite{Gel}, Huebschmann
\cite{Hue} and Fresse \cite{Fre} with the relation to Poisson
homology. The relations between canonical homology and Poisson
geometry were discussed by Vaisman in \cite{Vais}. The relation
between Poisson algebras and Hochschild homology of enveloping
algebras was investigated by Kassel \cite{Kass}.

\vspace{3mm}
\paragraph{\bf 2. The spectral sequence.}
All rings below are unital and all (bi)modules are unitary. Let
$R$ be a noetherian commutative ring of characteristic 0.
Unadorned tensor products mean tensor products over $R$. For an
increasing filtration ${\rm F}$, ${\rm F}_{p-1}\subset {\rm
F}_{p}$, we denote its associated gradation ${\rm Gr}$, ${\rm
Gr}=\bigoplus {\rm Gr}_{p}$, ${\rm Gr}_{p}={\rm F}_{p}/{\rm
F}_{p-1}$. For any graded abelian group $G$  we denote by $G_{p}$
its $p$-th homogeneous part.

\vspace{3mm}
\paragraph{\textbf{Definition 1.}}
A $\mathbb{Z}$-filtered associative $R$-algebra $A$ (resp. a
$\mathbb{Z}$-filtered bimodule over $A$, symmetric over $R$),
${\rm F}_{p-1}A\subset {\rm F}_{p}A$,
$\bigcap_{p\in\mathbb{Z}}{\rm F}_{p}A=0$,
$\bigcup_{p\in\mathbb{Z}}{\rm F}_{p}A=A$ (resp. ${\rm
F}_{p-1}M\subset {\rm F}_{p}M$, $\bigcap_{p\in\mathbb{Z}}{\rm
F}_{p}M=0$, $\bigcup_{p\in\mathbb{Z}}{\rm F}_{p}M=M$), is called
\textbf{almost commutative} (resp. \textbf{almost symmetric}) if
its associated graded algebra ${\rm Gr}A$ (resp. associated graded
bimodule ${\rm Gr}M$) is commutative (resp. symmetric). This means
that
\begin{align}
{\rm F}_{p_{0}}A\cdot{\rm F}_{p_{1}}A\subset {\rm
F}_{p_{0}+p_{1}}A,\ \ [ {\rm F}_{p_{0}}A, {\rm F}_{p_{1}}A]
\subset {\rm F}_{p_{0}+p_{1}-1}A
\end{align}
\begin{align}{\rm (resp.}\ {\rm F}_{p_{0}} A\cdot{\rm F}_{p_{1}}M, {\rm F}_{p_{0}}M\cdot{\rm F}_{p_{1}}A
\subset {\rm F}_{p_{0}+p_{1}}A,\ \ [ {\rm F}_{p_{0}}A, {\rm
F}_{p_{1}}M] \subset {\rm F}_{p_{0}+p_{1}-1}M).
\end{align}
\vspace{3mm}

On the Hochschild complex ${\rm C}_{\bullet}(A, M)=\bigoplus
_{k}{\rm C}_{k}(A, M)$, ${\rm C}_{k}(A, M)=M\otimes A^{\otimes k}$
we have an increasing filtration ${\rm F}_{p}$
\begin{align}
{\rm F}_{p}{\rm C}_{k}(A, M)=\sum_{p_{0}+\cdots+p_{k}\leq p} {\rm
F}_{p_{0}}M\otimes {\rm F}_{p_{1}}A\otimes\cdots\otimes {\rm
F}_{p_{k}}A.
\end{align}
This gives rise to a spectral sequence with ${\rm
E}^{1}_{p,q}={\rm H}_{p+q}({\rm Gr}A, {\rm Gr}M)_{p}$, the
homogeneous part of degree $p$ of the Hochschild homology ${\rm
H}_{p+q}({\rm Gr}A, {\rm Gr}M)$, converging to the Hochschild
homology ${\rm H}_{p+q}(A, M)$. For $M=A$ we obtain the spectral
sequence of Brylinski as in \cite{Bry}.

\vspace{3mm}
\paragraph{\bf 3. Canonical homology.}

\vspace{3mm}
\paragraph{\bf Definition 2.}
A   commutative graded algebra $B=\bigoplus_{p\in
\mathbb{Z}}B_{p}$, $B_{p}=0$ for $p\ll 0$ (resp. a symmetric
graded bimodule $N$ over a commutative graded algebra $B$,
$N=\bigoplus_{p\in \mathbb{Z}}N_{p}$, $N_{p}=0$ for $p\ll 0$), is
called a \textbf{ Poisson (graded) algebra} (resp. \textbf{Poisson
(graded) module} over a (graded) Poisson algebra $B$) if there is
given a Lie algebra structure on $B$
$$\{ -, -\} : B\otimes B\rightarrow B$$
(resp. a structure of a right module structure over the Lie
algebra $(B, \{ -,-\}) $
$$\{ -, -\} : N\otimes B\rightarrow N),$$
(with
$$\{B_{p_{0}}, B_{p_{1}}\} \subset B_{p_{0}+p_{1}-1}\ \ \ {\rm (resp.}\ \ \
\{N_{p_{0}}, B_{p_{1}}\} \subset N_{p_{0}+p_{1}-1}),$$ if they are
Poisson graded) such that for all $b_{0}, b_{1}, b_{2}\in B$
$$\{ b_{0}, b_{1}b_{2}\}=\{ b_{0}, b_{1}\}b_{2}+b_{1}\{ b_{0}, b_{2}\}$$
(resp. for all $n\in N, b_{1}, b_{2}\in B$
$$\{ nb_{1}, b_{2}\}=\{ n, b_{2}\}b_{1}+n\{ b_{1}, b_{2}\},\ \ \ \{n, b_{1}b_{2}\}=\{ n, b_{1}\}b_{2}+b_{1}\{ n, b_{2}\}.)$$

\vspace{3mm}
\paragraph{\bf Definition 3.} Let $N$ be a Poisson module over a Poisson algebra $B$. On the graded $R$-module
${\rm C}^{can}_{\bullet}(B, N)=\bigoplus _{k}C^{can}_{k}(B, N)$,
$C^{can}_{k}(B, N)=N\otimes_{B} \Omega^{k}_{B/R}$ one defines
\cite{Fre} the chain complex structure as follows:
$$\partial: C^{can}_{k}(B,
N)\rightarrow C^{can}_{k-1}(B, N),$$
$$\partial(n\otimes_{B} {\rm d}b_{1}\cdots{\rm d}b_{k})=$$
$$=\sum_{i=1}^{k}(-1)^{i-1}\{ n, b_{i}\}\otimes_{B}{\rm d}b_{1}\cdots\widehat{{\rm d}b_{i}}\cdots{\rm d}b_{k})$$
$$+\sum_{1\leq i, j\leq k}(-1)^{i+j}n\otimes_{B}{\rm d}\{ b_{i}, b_{j}\}{\rm d}b_{1}\cdots\widehat{{\rm d}b_{i}}\cdots\widehat{{\rm d}b_{j}}\cdots{\rm d}b_{k}).$$
One verifies that the boundary operator $\partial$ is well defined
and $\partial^{2}=0$. The homology ${\rm H}^{can}_{\bullet}(B, N)$
of this complex is called \textbf{canonical homology} of the
Poisson module $N$ over a Poisson algebra $B$. Note that if $B$
and $N$ are graded Poisson then $\partial$ is homogeneous of
degree (-1). Therefore in the Poisson graded case $k$-th canonical
chain and homology groups are graded in a canonical way.

\vspace{3mm}
\paragraph{\bf 4. The Hochschild-Kostant-Rosenberg isomorphism.}
We will use the simple observation \cite{Fre} that the
Hochschild-Kostant-Rosenberg isomorphism (see \cite{Bry},
\cite{Hoch}, \cite{Lod}) holds also in the case of coefficients in
a symmetric bimodule $N$ over a smooth commutative algebra $B$.
This means that the map
$$\beta: {\rm H}_{k}(B, N)\rightarrow N\otimes_{B}\Omega^{k}_{B/R},$$
\begin{align}
\beta(n\otimes b_{1}\otimes\cdots\otimes b_{k})=\frac{1}{k!}\
n\otimes_{B} db_{1}\cdots db_{k}
\end{align}
is an isomorphism, with the inverse $\gamma$:
\begin{align}
\gamma(n\otimes_{B} db_{1}\cdots db_{k})=\left[
\sum_{\sigma\in\mathfrak{S}_{k}}{\rm sgn}(\sigma)\ n\otimes
b_{\sigma(1)}\otimes\cdots\otimes b_{\sigma(k)}\right],
\end{align}
where the square bracket denotes the Hochschild homology class of
a cycle.

\vspace{3mm}
\paragraph{\bf 5. Hochschild and canonical homology.}
Given an almost commutative algebra $A$ (resp. an almost symmetric
bimodule $M$ over $A$), $B={\rm Gr}A$ (resp. $N={\rm Gr}M$) has a
canonical structure of a graded Poisson algebra (resp. a graded
Poisson module over $B$) with the Poisson structure  defined as
follows: for $b_{0}=a_{0}+{\rm F}_{p_{0}-1}A$, $a_{0}\in {\rm
F}_{p_{0}}A$, $b_{1}=a_{1}+{\rm F}_{p_{1}-1}A$, $a_{1}\in {\rm
F}_{p_{1}}A$, $n=m+{\rm F}_{p_{0}-1}M$, $m\in {\rm F}_{p_{0}}M$
\begin{align*}
\{ b_{0}, b_{1}\} & =[a_{0}, a_{1}]+{\rm F}_{p_{0}+p_{1}-2}A,\\
\{ n, b_{1}\} & =[m, a_{1}]+{\rm F}_{p_{0}+p_{1}-2}M.
\end{align*}

Moreover we have the following theorem generalizing Theorem 3.1.1
of \cite{Bry}. In the proof we follow the lines of the beautiful
proof of \cite{Bry}, improving a little misprint in the original
proof ($a_{i}$ instead of $a_{0}$ in the formula (II) in
\cite{Bry}). At first sight this (spoiled) structure of the
formula (II) makes our generalization impossible, but after this
minor correction everything can be adapted verbatim.

\begin{theorem} Assume that the above filtrations are bounded and exhaustive, and $B={\rm Gr}A$ is smooth over $R$.
Then for any $q\geq 0$ the Hochschild-Kostant-Rosenberg
isomorphism induces an isomorphism of complexes
\begin{align}
\beta : ({\rm E}^{1}_{p,q}(A, M), d^{1}_{p,q})\rightarrow({\rm
C}^{can}_{p+q}(B, N)_{p},
\partial ),
\end{align}
where $d^{1}_{p,q}: {\rm E}^{1}_{p,q}(A, M)\rightarrow {\rm
E}^{1}_{p-1,q}(A, M)$ is the differential in the spectral
sequence.

In particular
\begin{align}
{\rm E}^{2}_{p,q}(A, M)\cong {\rm H}^{can}_{p+q}(B, N)_{p}.
\end{align}
\end{theorem}

\textit{Proof.} It is enough to prove that $\beta\circ
d^{1}\circ\gamma=\partial$. Now, the $R$-module $C^{can}_{k}(B,
N)_{p}=(N\otimes_{B} \Omega^{k}_{B/R})_{p}$ is generated by
elements of the form $n\otimes_{B} {\rm d}b_{1}\cdots{\rm
d}b_{k}$, where $n=m+{\rm F}_{p_{0}-1}M$, $m\in {\rm F}_{p_{0}}M$,
$b_{i}=a_{i}+{\rm F}_{p_{i}-1}A$, $a_{i}\in {\rm F}_{p_{i}}A$,
$p_{0}+\cdots+p_{k}=p$. First, we have
$$\gamma(n\otimes_{B} db_{1}\cdots db_{k})=\left[
\sum_{\sigma\in\mathfrak{S}_{k}}{\rm sgn}(\sigma)\ n\otimes
b_{\sigma(1)}\otimes\cdots\otimes b_{\sigma(k)}\right].$$ The
cycle on the right hand side lives in ${\rm C}_{k}(B, N)_{p}$ and
lifts to the chain
\begin{align}
\sum_{\sigma\in\mathfrak{S}_{k}}{\rm sgn}(\sigma)\ m\otimes
a_{\sigma(1)}\otimes\cdots\otimes a_{\sigma(k)}\in {\rm F}_{p}{\rm
C}_{k}(A, M).\label{chain}\end{align} Its Hochschild boundary is
the sum of three terms (I), (II), (III), with
\begin{align*}
{\rm (I)} & =\sum_{\sigma\in\mathfrak{S}_{k}}{\rm sgn}(\sigma)\
ma_{\sigma(1)}\otimes a_{\sigma(2)}\otimes\cdots\otimes a_{\sigma(k)},\\
{\rm (II)} & =\sum_{\sigma\in\mathfrak{S}_{k}}\sum_{1\leq i<k}{\rm
sgn}(\sigma)(-1)^{i}
m\otimes a_{\sigma(1)}\otimes\cdots \otimes a_{\sigma(i)}a_{\sigma(i+1)}\otimes\cdots \otimes a_{\sigma(k)},\\
{\rm (III)} & =\sum_{\sigma\in\mathfrak{S}_{k}}{\rm
sgn}(\sigma)(-1)^{k} a_{\sigma(1)}m\otimes
a_{\sigma(1)}\otimes\cdots\otimes a_{\sigma(k-1)}.
\end{align*}
Since the chain (\ref{chain}) is a lift of a Hochschild cycle in
${\rm C}_{k}(B, N)_{p}$ its Hochschild boundary lives in ${\rm
F}_{p-1}{\rm C}_{k-1}(A, M)$. Now we are to compute the image of
this Hochschild boundary in ${\rm Gr}_{p-1}{\rm C}_{k-1}(A,
M)={\rm C}_{k-1}(B, N)_{p-1}$.

First, transforming $\sigma\in\mathfrak{S}_{k}$ to $\sigma\tau$,
where $\tau$ is a cyclic permutation, we can rewrite (I) as
follows
\begin{align*}
{\rm (I)} & =\sum_{\sigma\in\mathfrak{S}_{k}}(-1)^{k+1}{\rm
sgn}(\sigma)\ ma_{\sigma(k)}\otimes
a_{\sigma(1)}\otimes\cdots\otimes a_{\sigma(k-1)}.
\end{align*}
Since
\begin{align*}
{\rm (I)}+{\rm (III)} &
=\sum_{\sigma\in\mathfrak{S}_{k}}(-1)^{k+1}{\rm sgn}(\sigma)\
[m,a_{\sigma(k)}]\otimes a_{\sigma(1)}\otimes\cdots\otimes
a_{\sigma(k-1)}
\end{align*}
we have ${\rm (I)}+{\rm (III)}\in {\rm F}_{p-1}{\rm C}_{k-1}(A,
M)$; its image in ${\rm Gr}_{p-1}{\rm C}_{k-1}(A, M)={\rm
C}_{k-1}(B, N)_{p-1}$ is equal to
\begin{align}
\sum_{\sigma\in\mathfrak{S}_{k}}(-1)^{k+1}{\rm sgn}(\sigma)\ \{
n,b_{\sigma(k)}\}\otimes b_{\sigma(1)}\otimes\cdots\otimes
b_{\sigma(k-1)}.\label{I+III}
\end{align}
Second, transforming  $\sigma\in\mathfrak{S}_{k}$ to $\sigma
s_{h}$, where $s_{h}$ is a transposition which exchanges $h$ and
$(h+1)$, we can rewrite (II) as follows
\begin{align*}
{\rm (II)} &
=\frac{1}{2}\sum_{\sigma\in\mathfrak{S}_{k}}\sum_{1\leq h<k}{\rm
sgn}(\sigma)(-1)^{h} m\otimes a_{\sigma(1)}\otimes\cdots \otimes
[a_{\sigma(h)},a_{\sigma(h+1)}]\otimes\cdots \otimes a_{\sigma(k)}
\end{align*}
which also lives in ${\rm F}_{p-1}{\rm C}_{k-1}(A, M)$; its image
in ${\rm Gr}_{p-1}{\rm C}_{k-1}(A, M)={\rm C}_{k-1}(B, N)_{p-1}$
is equal to
\begin{align}
\frac{1}{2}\sum_{\sigma\in\mathfrak{S}_{k}}\sum_{1\leq h<k}{\rm
sgn}(\sigma)(-1)^{h} n\otimes b_{\sigma(1)}\otimes\cdots \otimes
\{ b_{\sigma(h)},b_{\sigma(h+1)}\} \otimes\cdots \otimes
b_{\sigma(k)}.\label{II}
\end{align}
Adding (\ref{I+III}) and (\ref{II}) we obtain
$(d^{1}\circ\gamma)(n\otimes_{B} db_{1}\cdots db_{k})$. It remains
to apply $\beta$ to this.

Now, for the image of the sum (\ref{I+III}) under $\beta$ notice
that all $\sigma$'s with $\sigma(k)=i$ fixed, give the same value
for $\beta(\{ n,b_{\sigma(k)}\}\otimes
b_{\sigma(1)}\otimes\cdots\otimes b_{\sigma(k-1)})$ equal to
\begin{align*}
\frac{1}{(k-1)!}(-1)^{k-i}\{ n,b_{i}\}\otimes_{B} {\rm d
}b_{1}\cdots \widehat{{\rm d }b_{i}} \cdots {\rm d }
b_{\sigma(k)}.
\end{align*}
Since there are $(k-1)!$ such permutations, substituting this
value, common for each $\sigma$ with $\sigma(k)=i$, to the image
of (\ref{I+III}) we obtain
\begin{align}
\sum_{1\leq i\leq k}(-1)^{i-1}\{ n,b_{i}\}\otimes_{B} {\rm d
}b_{1}\cdots \widehat{{\rm d }b_{i}} \cdots {\rm d }
b_{\sigma(k)}. \label{A}\end{align}

Next, for the image of (\ref{II}) under $\beta$ notice that all
pairs $(\sigma , h)$ with the set $\{ \sigma(h),\sigma(h+1)\}$
equal to a fixed set $\{ i,j\}$ (say $i<j$), give the same value
for $\beta({\rm sgn}(\sigma)(-1)^{h} n\otimes
b_{\sigma(1)}\cdots\otimes \{ b_{\sigma(h)}, b_{\sigma(h+1)}
\}\otimes\cdots\otimes b_{\sigma(k)})$ equal to
\begin{align*}
\frac{1}{(k-1)!}(-1)^{i+j} n\otimes_{B}{\rm d}\{ b_{i},b_{j}
\}{\rm d }b_{1}\cdots \widehat{{\rm d }b_{i}} \cdots \widehat{{\rm
d }b_{j}}\cdots{\rm d } b_{k}.
\end{align*}
Since there are $2(k-1)!$ such pairs, substituting this value,
common for each $(\sigma, h)$ with the set $\sigma(h)=i$ equal to
a fixed set, to the image of (\ref{II}) we obtain
\begin{align}
\sum_{1\leq i,j\leq k}(-1)^{i+j}n\otimes_{B}{\rm d}\{ b_{i},b_{j}
\}{\rm d }b_{1}\cdots \widehat{{\rm d }b_{i}} \cdots \widehat{{\rm
d }b_{j}}\cdots{\rm d } b_{k}. \label{B}\end{align} Taking the sum
of (\ref{A}) and (\ref{B}) we obtain
\begin{align*}
(\beta\circ d^{1}\circ\gamma)(n\otimes_{B} db_{1}\cdots
db_{k})=\partial(n\otimes_{B} db_{1}\cdots db_{k}).\  \Box
\end{align*}

\vspace{3mm}
\paragraph{\bf 6. Maximal commutative subalgebras and almost commutative algebras.}

\vspace{3mm}
\paragraph{\bf Definition 4.} Let $C$ be a
maximal commutative subalgebra of $A$ and $M$ be a bimodule over
$A$. For $c\in C$, $m\in M$ we define an operation ${\rm
ad}_{c}(m):=[c,m]$ and an increasing $\mathbb{N}$-filtration on
$M$
\begin{align}
{\rm F}^{C}_{p}M:=\{m\in M \mid \forall_{c\in C}\ \  {\rm
ad}^{p+1}_{c}(m)= 0\}.\label{filt}
\end{align}
Since for all $c, c'\in C$ $[{\rm ad}_{c}, {\rm ad}_{c'}]={\rm
ad}_{[c,c']}=0,$ the multilinear map
\begin{align}
(c_{0},\ldots,c_{p})\mapsto{\rm ad}_{c_{0}}\ldots {\rm
ad}_{c_{p}}(m)\label{polar}
\end{align}
is symmetric in $(c_{0},\ldots,c_{p})$. Therefore (\ref{polar}) is
a symmetric $(p+1)$-linear form in $(c_{0},\ldots,c_{p})$
corresponding, via the polarization formula, to the homogeneous
polynomial of degree $(p+1)$ in $c$
\begin{align}
c\mapsto{\rm ad}_{c}^{p+1}(m).\label{polyn}
\end{align}
Therefore the filtration ${\rm F}^{C}$ can be rewritten
equivalently as follows
\begin{align}
{\rm F}^{C}_{p}M:=\{m\in M \mid \forall_{c_{0},\ldots,c_{p}\in C}\
\  {\rm ad}_{c_{0}}\ldots {\rm ad}_{c_{p}}(m)= 0\}.\label{filt'}
\end{align}
We call ${\rm D}^{C}A:=\bigcup_{p}{\rm F}^{C}_{p}A$ (resp. ${\rm
D}^{C}M:=\bigcup_{p}{\rm F}^{C}_{p}M$) the \textbf{differential
hull} of $C$ in $A$ (resp. of the centralizer of $C$ in $M$). Note
that ${\rm D}^{C}M$ inherits the filtration from $M$ and has the
same associated gradation.

The following example justifies the name of the differential hull.

\vspace{3mm}
\paragraph{\bf Example 1.} Let us take two left $A$-modules $P$ and
$Q$, and form the $A$-bimodule $M:={\rm Hom}_{R}(P, Q)$. Then by
(\ref{filt'}) we have
$${\rm F}^{C}_{p}M={\rm Diff}_{p}^{C/R}(P, Q),\ \ {\rm Gr}^{C}_{p}M={\rm Smbl}_{p}^{C/R}(P, Q),$$
where ${\rm Diff}_{p}^{C/R}(P, Q)$ (resp. ${\rm Smbl}_{p}^{C/R}(P,
Q)$) denotes the $C$-bimodule of differential operators of order
$p$ from $P$ to $Q$ (resp. the symmetric $C$-bimodule of their
principal symbols).

\begin{theorem}
Given a maximal commutative subalgebra $C$ of $A$, the above
filtration makes ${\rm D}^{C}A$ almost commutative and a bimodule
${\rm D}^{C}M$ over ${\rm D}^{C}A$ almost symmetric.
\end{theorem}

\textit{Proof.}
 Let $m\in {\rm F}^{C}_{p_{0}}M$, $a\in
{\rm F}^{C}_{p_{1}}A$.

On the right hand side of the identities
\begin{align*}
\frac{1}{(p_{0}+p_{1}+1)!}{\rm
ad}^{p_{0}+p_{1}+1}_{c}(ma)=\sum_{i+j=p_{0}+p_{1}+1}\frac{1}{i!}{\rm
ad}^{i}_{c}(m)\frac{1}{j!}{\rm ad}^{j}_{c}(a)
\end{align*}
\begin{align*}
\frac{1}{(p_{0}+p_{1}+1)!}{\rm
ad}^{p_{0}+p_{1}+1}_{c}(am)=\sum_{i+j=p_{0}+p_{1}+1}\frac{1}{j!}{\rm
ad}^{j}_{c}(a)\frac{1}{i!}{\rm ad}^{i}_{c}(m)
\end{align*}
at least one of the two factors in every summand is zero, since
either $i\geq p_{0}+1$ or $j\geq p_{1}+1$. This proves that
\begin{align}
{\rm F}^{C}_{p_{0}}M\cdot{\rm F}^{C}_{p_{1}}A,\ {\rm
F}^{C}_{p_{1}}A\cdot{\rm F}^{C}_{p_{0}}M\subset {\rm
F}^{C}_{p_{0}+p_{1}}M.\label{mult}
\end{align}

Observe now that ${\rm F}^{C}_{0}A=C$, since $C$ is a maximal
commutative subalgebra in $A$, and $[{\rm F}^{C}_{0}M, C]=0$.
Next, if $c\in C$, $m\in {\rm F}^{C}_{p}M$  then
\begin{align*}
{\rm ad}^{p+1-i}_{c}({\rm ad}^{i}_{c}(m))={\rm ad}^{p+1}_{c}(m)=0
\end{align*}
which means that
\begin{align}
{\rm ad}^{i}_{c}({\rm F}^{C}_{p}M)\subset {\rm
F}^{C}_{p-i}M.\label{decr}
\end{align}

On the right hand side of the identity
\begin{align*}
\frac{1}{(p_{0}+p_{1})!}{\rm
ad}^{p_{0}+p_{1}+1}_{c}([m,a])=\sum_{i+j=p_{0}+p_{1}}\left[
\frac{1}{i!}{\rm ad}^{i}_{c}(m), \frac{1}{j!}{\rm
ad}^{j}_{c}(a)\right]
\end{align*}
all summands are zero, since the following implications hold:
\begin{eqnarray*}
i<p_{0} & \ \Rightarrow\ & j>p_{1}\ \Rightarrow\ {\rm
ad}^{j}_{c}(a)=0,\\
i=p_{0} & \ \Rightarrow\ & j=p_{1}\
\stackrel{(\ref{decr})}{\Rightarrow}\ {\rm ad}^{i}_{c}(m)\in {\rm
F}^{C}_{0}M,\  {\rm
ad}^{j}_{c}(a)\in {\rm F}^{C}_{0}A=C,\\
i>p_{0} & \ \Rightarrow\ & {\rm ad}^{i}_{c}(m)=0.
\end{eqnarray*}
This proves that
\begin{align}
\left[ {\rm F}^{C}_{p_{0}}M, {\rm F}^{C}_{p_{1}}A,\right] \subset
{\rm F}^{C}_{p_{0}+p_{1}-1}M.\label{comm}
\end{align}
Taking $M=A$ in (\ref{mult}) and (\ref{comm}) we see that the
filtration ${\rm F}^{C}$ makes ${\rm D}^{C}A$ almost commutative
and, also by (\ref{mult}) and (\ref{comm}), ${\rm D}^{C}M$ almost
symmetric over ${\rm D}^{C}A$. $\Box$

\vspace{3mm}
\paragraph{\bf 7. Poisson geometry of the differential hull.}
To discuss the Poisson geometry arising from the differential hull
we need to generalize the notion of an involutive distribution and
a sheaf with a flat connection along an involutive distribution.

Let $C$ be a smooth commutative $R$-algebra and let
$\Theta^{C/R}={\rm Hom}_{C}(\Omega^{1}_{C/R}, C)$ denote the
relative tangent module of $C$ over $R$.

Every $f\in {\rm Hom}_{C}({\rm Sym}^{p}_{C}\Omega^{1}_{C/R},
C)={\rm Sym}^{p}_{C}\Theta^{C/R}$ can be regarded as a symmetric
$R$-linear $C$-valued $p$-form on $C$, which is a derivation with
respect to every linear argument. Using this fact we can define
the canonical graded Poisson algebra structure of the $R$-algebra
${\rm Sym}_{C}\Theta^{C/R}=\bigoplus_{p\geq 0} {\rm
Sym}^{p}_{C}\Theta^{C/R}$ of polynomial functions on the relative
cotangent bundle of the scheme ${\rm Spec}C$ over $R$ as follows:
for $f_{i}\in {\rm Hom}_{C}({\rm Sym}^{p_{i}}_{C}\Omega^{1}_{C/R},
C)$, $i=0, 1$,
\begin{align}
(f_{0}f_{1})(c_{1},\ldots,c_{p_{0}+p_{1}})
:=\label{pro}\end{align}
\begin{align}
  \frac{1}{(p_{0}+p_{1})!} \sum_{\sigma\in\mathfrak{S}_{p_{0}+p_{1}}} & f_{0}(c_{\sigma(1)},\ldots,c_{\sigma(p_{0})})f_{1}(c_{\sigma(p_{0}+1)},\ldots,c_{\sigma(p_{0}+p_{1})}),\nonumber
\end{align}
\begin{align}
\{ f_{0}, f_{1}\}(c_{1},\ldots,c_{p_{0}+p_{1}-1})
:=\label{bra}\end{align}
\begin{align}
  \frac{1}{(p_{0}+p_{1}-1)!}  \sum_{\sigma\in\mathfrak{S}_{p_{0}+p_{1}-1}}
  ( &
  p_{0}f_{0}(c_{\sigma(1)},\ldots,c_{\sigma(p_{0}-1)},f_{1}(c_{\sigma(p_{0})},\ldots,c_{\sigma(p_{0}+p_{1}-1)}))\nonumber\\
 - & p_{1}f_{1}(c_{\sigma(1)},\ldots,c_{\sigma(p_{1}-1)},f_{0}(c_{\sigma(p_{1})},\ldots,c_{\sigma(p_{0}+p_{1}-1)}))).
\end{align}

\vspace{3mm}
\paragraph{\bf Definition 5.} A \textbf{\textit{nonlinear involutive distribution}} on
a scheme ${\rm Spec}C$ over $R$ is a graded Poisson subalgebra $B$
of the algebra ${\rm Sym}_{C}\Theta^{C/R}$ of polynomial functions
on the relative cotangent bundle of the scheme ${\rm Spec}C$ over
$R$ such that $B_{0}=C$. A \textbf{\textit{graded sheaf with a
flat connection}} along the nonlinear involutive distribution $B$
is a graded Poisson module $N$ over $B$ together with a structure
of a graded Poisson module over $B$ on $N_{0}\otimes_{C}{\rm
Sym}_{C}\Theta^{C/R}$ and an embedding of $N$ as a graded Poisson
submodule of $N_{0}\otimes_{C}{\rm Sym}_{C}\Theta^{C/R}$.

\vspace{3mm}
\paragraph{\bf Example 2.} Let $L\subset \Theta^{C/R}={\rm Der}_{R}(C,C)$ be a finitely
generated projective $C$-submodule which is also an
$R$-Lie-subalgebra. It describes an involutive distribution on
${\rm Spec}C$ over $R$. The $C$-linear embedding $L\hookrightarrow
\Theta^{C/R}$ is equivalent to the grading preserving embedding of
graded $C$-algebras
\begin{align}
{\rm Sym}_{C}L\hookrightarrow{\rm Sym}_{C}\Theta^{C/R}.\label{lin}
\end{align}
The Lie subalgebra structure defines on the image of (\ref{lin}) a
structure of a graded Poisson subalgebra uniquely determined by
the following brackets for $\theta\in \Theta^{C/R}$, $c\in C$,
$l\in L$
\begin{align} \{  \theta,
c\}=  \theta(c),\ \ \{  \theta,  l\}=[\theta, l].\label{bra'}
\end{align}
This shows that an involutive distribution is an instance of a
nonlinear involutive distribution. On the other hand, every
nonlinear involutive distribution of the form ${\rm Sym}_{C}L$ for
some finitely generated projective $C$-module $L$ defines on $L$ a
structure of an involutive distribution. Assume now that $N_{0}$
is a $C$-module equipped with a flat connection along an
involutive distribution $L\subset \Theta^{C/R}$
\begin{align}\nabla:N_{0}\otimes_{R}L & \rightarrow N,\ \
n_{0}\otimes l \mapsto \nabla_{l}n_{0},\label{con}
\end{align}
\begin{align}
\nabla_{cl}n_{0}=c\nabla_{l}n_{0},\ \
\nabla_{l}(n_{0}c)=(\nabla_{l}n_{0})c+n_{0}l(c),\ \
[\nabla_{l'},\nabla_{l}]=\nabla_{[l,l']}.\nonumber
\end{align}
Then the brackets
\begin{align}\{ n_{0}\otimes 1, 1\otimes c\}=0,\ \ \{ n_{0}\otimes 1,1\otimes
l\}=-\nabla_{l}n_{0}\otimes 1\label{bra''}
\end{align}
determine uniquely structures of  graded Poisson modules over
${\rm Sym}_{C}L$ on $N:=N_{0}\otimes_{C}{\rm Sym}_{C}L$ and
$N_{0}\otimes_{C}{\rm Sym}_{C}\Theta^{C/R}$, and an embedding (of
graded Poisson modules over ${\rm Sym}_{C}L$) $N\hookrightarrow
N_{0}\otimes_{C}{\rm Sym}_{C}\Theta^{C/R}$. In this way $N$
becomes a graded sheaf with a flat connection along the nonlinear
involutive distribution ${\rm Sym}_{C}L$.
\begin{theorem} Let $C$ be a smooth commutative $R$-subalgebra of
an associative $R$-algebra $A$ and let $M$ be an arbitrary
$R$-symmetric $A$-bimodule. Then ${\rm Gr}^{C}A$ is a nonlinear
involutive distribution on a scheme ${\rm Spec}C$ over $R$ and
${\rm Gr}^{C}M$ is a graded sheaf with a flat connection along
${\rm Gr}^{C}A$.
\end{theorem}

\textit{Proof.} Since an $A$-bimodule $M$ is symmetric as an
$R$-bimodule and $C$ is a commutative $R$-subalgebra in $A$ $M$ is
a left module over $C\otimes_{R}C$, where $(c_{0}\otimes
c_{1})m:=c_{0}m c_{1}$. Let us consider the kernel $I$ of the
multiplication map $C\otimes_{R}C\rightarrow C$. This is an ideal
in $C\otimes_{R}C$ generated by elements of the form $c\otimes
1-1\otimes c$. Since $(c\otimes 1-1\otimes c)m={\rm ad}_{c}(m)$
the filtration (\ref{filt'}) gives rise, for any $k\leq p$, to the
embedding
\begin{align}{\rm
Gr}^{C}_{p}M  \hookrightarrow {\rm Hom}_{C}(I^{k}/I^{k+1}, {\rm
Gr}^{C}_{p-k}M),\label{emb}\end{align}
$$m+{\rm F}^{C}_{p-1}M \mapsto $$
$$((c_{1}\otimes 1-1\otimes c_{1})\cdots (c_{k}\otimes
1-1\otimes c_{k})+I^{k+1}\mapsto \frac{(-1)^{k}}{k!}{\rm
ad}_{c_{1}}\ldots {\rm ad}_{c_{k}}(m)+{\rm F}^{C}_{p-k-1}M).$$
Since for $C$ smooth over $R$ one has the isomorphism of symmetric
$C$-bimodules
\begin{align}
{\rm Sym}^{k}_{C}\Omega^{1}_{C/R}\stackrel{\cong}{\rightarrow}
I^{k}/I^{k+1}\label{sym}
\end{align}
and ${\rm Gr}^{C}_{0}M$ is a symmetric $C$-bimodule the latter
embedding can be rewritten as
\begin{align} {\rm
Gr}^{C}_{p}M & \hookrightarrow {\rm Hom}_{C}({\rm
Sym}^{k}_{C}\Omega^{1}_{C/R}, {\rm Gr}^{C}_{p-k}M)={\rm
Gr}^{C}_{p-k}M\otimes_{C}{\rm Sym}^{k}_{C}\Theta^{C/R}.\label{mod}
\end{align}

 In particular, since ${\rm Gr}^{C}_{0}A=C$, we obtain the
embedding
\begin{align}
{\rm Gr}^{C}_{p}A & \hookrightarrow {\rm Hom}_{C}({\rm
Sym}^{p}_{C}\Omega^{1}_{C/R}, C)={\rm
Sym}^{p}_{C}\Theta^{C/R}.\label{fun}
\end{align}

Then the embeddings (\ref{fun}) define a grading preserving
embedding of graded Poisson algebras
\begin{align}
{\rm Gr}^{C}A & \hookrightarrow {\rm
Sym}_{C}\Theta^{C/R}.\label{fun'}
\end{align}

The canonical structure of a graded module over the graded algebra
${\rm Sym}_{C}\Theta^{C/R}$ on ${\rm Gr}_{0}^{C}M\otimes_{C}{\rm
Sym}_{C}\Theta^{C/R}$ is consistent through (\ref{fun'}) with the
following structure of a graded Poisson module over ${\rm
Gr}^{C}A$ on ${\rm Gr}_{0}^{C}M\otimes_{C}{\rm
Sym}_{C}\Theta^{C/R}$ identified with $\bigoplus_{p\geq 0} {\rm
Hom}_{C}({\rm Sym}^{p}_{C}\Omega^{1}_{C/R}, {\rm Gr}_{0}^{C}M)$:
for all $s_{0}\in {\rm Hom}_{C}({\rm
Sym}^{p_{0}}_{C}\Omega^{1}_{C/R}, {\rm Gr}_{0}^{C}M))$ and
$f_{1}\in {\rm Hom}_{C}({\rm Sym}^{p_{1}}_{C}\Omega^{1}_{C/R}, C)$
\begin{align}
(s_{0}f_{1})(c_{1},\ldots,c_{p_{0}+p_{1}})
:=\label{pro'}\end{align}
\begin{align}
  \frac{1}{(p_{0}+p_{1})!} \sum_{\sigma\in\mathfrak{S}_{p_{0}+p_{1}}} & s_{0}(c_{\sigma(1)},\ldots,c_{\sigma(p_{0})})f_{1}(c_{\sigma(p_{0}+1)},\ldots,c_{\sigma(p_{0}+p_{1})}),\nonumber
\end{align}
\begin{align}
\{ s_{0}, b_{1}\}(c_{1},\ldots,c_{p_{0}+p_{1}-1})
:=\label{bra'''}\end{align}
\begin{align}
  \frac{1}{(p_{0}+p_{1}-1)!}  \sum_{\sigma\in\mathfrak{S}_{p_{0}+p_{1}-1}}
  (  &
  p_{0}s_{0}(c_{\sigma(1)},\ldots,c_{\sigma(p_{0}-1)},b_{1}(c_{\sigma(p_{0})},\ldots,c_{\sigma(p_{0}+p_{1}-1)}))\nonumber\\
+ & p_{1}\{
s_{0}(c_{\sigma(p_{1})},\ldots,c_{\sigma(p_{0}+p_{1}-1)}),
b_{1}(c_{\sigma(1)},\ldots,c_{\sigma(p_{1}-1)},-)\}),\nonumber
 \end{align}
where $\{-,-\}$ on the right hand side is the bracket ${\rm
Gr}^{C}_{0}M\otimes_{R}{\rm Gr}_{1}^{C}A  \rightarrow{\rm
Gr}_{0}^{C}M$. Then the embeddings (\ref{mod}) for $k=p$ define an
embedding of graded Poisson modules over ${\rm Gr}^{C}A$
\begin{align} {\rm
Gr}^{C}M & \hookrightarrow {\rm Gr}^{C}_{0}M\otimes_{C}{\rm
Sym}_{C}\Theta^{C/R}.\label{mod'}
\end{align}
$\Box$

There are also interesting examples of non-smooth maximal
commutative subalgebras.

\vspace{3mm}
\paragraph{\bf Example 3.} Let $A=M_{n}(k)$ be the $(n\times n)$ matrix algebra
over an algebraically closed field $k$ of characteristic 0. Among
its many non-isomorphic maximal commutative subalgebras \cite{Sup}
one has the following two extremes:

1) Diagonal subalgebra $C\cong k^{n}=k\times\cdots\times k$. It is
smooth (semisimple) and the filtration stabilizes at $C$, hence
${\rm D}^{C}A=C$, ${\rm Gr}^{C}A=C$,

2) Subalgebra $C\cong k[x]/(x^{n+1})$ generated by the nilpotent
$(n\times n)$ Jordan block. Then the filtration is exhaustive,
${\rm D}^{C}A={\rm Diff}_{R}C$, the algebra of differential
operators of the commutative $k$-algebra $C$, and  ${\rm Gr}^{C}A$
is isomorphic as a graded Poisson algebra to the polynomial
algebra of the $n$th infinitesimal neighborhood of the nilpotent
cone in ${\rm sl}^{*}_{2}$ with its canonical
Kirillov-Kostant-Souriau Poisson structure and the grading such
that ${\rm deg}(e)=0$, ${\rm deg}(h)=1$, ${\rm deg}(f)=2$
\cite{Masz}.
\vspace{3mm}
\paragraph{\bf 8. Deformation of the differential hull to a nonlinear involutive distribution.}
According to the construction of Gerstenhaber \cite{Ger} we apply
a slight generalization of the Rees algebra to our almost
commutative algebra ${\rm D}^{C}A$ and an almost symmetric
bimodule ${\rm D}^{C}M$.

\vspace{3mm}
\paragraph{\bf Definition 6.} For any bimodule $M$ over an associative
$R$-algebra $A$, and a maximal commutative subalgebra $C$ of $A$,
we define
\begin{align}
{\rm Def}M:=\sum_{p}{\rm F}^{C}_{p}M\cdot t^{p}\subset
M[t]=M\otimes_{R}R[t]. \label{def}
\end{align}
By Theorem 3, according to \cite{Ger}, we have the following
corollary.
\begin{corollary} In the situation of Definition 6 the following
holds.

1) ${\rm Def}A$ is a subalgebra in $A[t]$ and ${\rm Def}M$ is a
sub-bimodule of a ${\rm Def}A$-bimodule $M[t]$,

2) At $t=1$ ${\rm Def}A$ (resp. ${\rm Def}M$) specializes to the
almost commutative algebra ${\rm D}^{C}A$ (resp. to the almost
symmetric ${\rm Def}A$-bimodule ${\rm Def}M$),

3) At $t=0$ ${\rm Def}A$ (resp. ${\rm Def}M$) specializes to the
nonlinear involutive distribution ${\rm Gr}^{C}A$  (resp. to a
graded sheaf with a flat connection along ${\rm Gr}^{C}A$) on the
scheme  ${\rm Spec}C$ over $R$.
\end{corollary}

\vspace{3mm}
\paragraph{\bf 9. Hochschild homology of the differential hull.}

 As a corollary of Theorem 1 and Theorem 2 we obtain
the following

\begin{theorem}
Let $C$ be a maximal commutative subalgebra of an associative
$R$-algebra $A$ and $M$ be a bimodule over $A$. Then

1) $C$ defines filtrations ${\rm F}^{C}$ making ${\rm D}^{C}A$
almost commutative and ${\rm D}^{C}M$ almost symmetric over ${\rm
D}^{C}A$, such that ${\rm Gr}^{C}A$ is a graded Poisson algebra
and ${\rm Gr}^{C}M$ a graded Poisson module over ${\rm Gr}^{C}A$.

2) These filtrations give rise to a spectral sequence converging
to the Hochschild homology
\begin{align*}
{\rm E}^{r}_{p,q}\ \Rightarrow {\rm H}_{p+q}({\rm D}^{C}A, {\rm
D}^{C}M).
\end{align*}

3) If ${\rm Gr}^{C}A$ is smooth over $R$ then
\begin{align*}
{\rm E}^{2}_{p,q}\cong {\rm H}^{can}_{p+q}({\rm Gr}^{C}A, {\rm
Gr}^{C}M)_{p}.
\end{align*}
\end{theorem}

\end{document}